\documentclass[12pt]{amsart} 
  \usepackage{latexsym} 
 \usepackage[all]{xy} 
  \usepackage{amsfonts} 
  \usepackage{amsthm} 
  \usepackage{amsmath} 
  \usepackage{amssymb} 

  \def\<{{\langle}} 
  \def\>{{\rangle}} 
   
  \def\eps{\varepsilon}

  \def\note#1{{}} 
  \def\can{{\rm can}} 
   
  \def\note#1{}

  \def\cC{{\mathfrak C}} 
   
  \def\cD{{\mathfrak D}}

  \def\cO{{\mathcal O}}

  \def\can{{\rm can}}

  \def\beq{\begin{equation}} 
  \def\eeq{\end{equation}} 
  \def\DC{{\Delta_\cC}} 
  \def \eC{{\eps_\cC}}

  \def\id{{I}}

  \def\ot{{\otimes}}

  \newcounter{zlist} 
  \newenvironment{zlist}{\begin{list}{(\arabic{zlist})}{ 
  \usecounter{zlist}\leftmargin2.5em\labelwidth2em\labelsep0.5em 
  \topsep0.6ex
  \parsep0.3ex plus0.2ex minus0.1ex}}{\end{list}}

  \newcounter{blist}

  \newcounter{rlist}

  \def\Label#1{\label{#1}\ifmmode\llap{[#1] }\else 
  \marginpar{\smash{\hbox{\tiny [#1]}}}\fi} 
  \def\Label{\label}

  \newtheorem{proposition}{Proposition}[section]

  \newtheorem{corollary}[proposition]{Corollary} 
  \newtheorem{theorem}[proposition]{Theorem} 

  \theoremstyle{definition} 
  \newtheorem{definition}[proposition]{Definition}

  \theoremstyle{remark}

  \newcounter{c} 
   
  \newcommand{\etyk}[1]{\vspace{-7.4mm}$$\begin{equation}\Label{#1} 
  \addtocounter{c}{1}} 
  \renewcommand{\]}{\ifnum \value{c}=1 $$\else \end{equation}\fi} 
\begin{document} 

  \title{Complete coverings and Galois corings} 
  \author{Tomasz Brzezi\'nski}    \author{Adam P.\ Wrightson} 
  \address{ Department of Mathematics, University of Wales Swansea, 
  Singleton Park, \newline\indent  Swansea SA2 8PP, U.K.} 
  \email{T.Brzezinski@swansea.ac.uk (T Brzezinski)} \email{maapw@swansea.ac.uk
  (AP Wrightson)}  
  \begin{abstract} 
 It is shown that any finite complete covering of a non-commutative  algebra in 
the sense of Calow and Matthes (J.\ Geom.\ Phys.\ 32 (2000), 114--165) gives rise to a 
 Galois coring.
   \end{abstract} 
  \maketitle 

  \section{Introduction}
  Early in the development of non-commutative geometry 
  it was realised that an 
algebraic structure underlying the geometric notion
  of a non-commutative 
principal bundle is provided by a Hopf-Galois extension 
  (cf.\ \cite{Sch:pri}). 
Motivated by the importance of principal bundles in
  physics, where they 
provide a geometric framework for gauge theories, the 
  foundations of the 
theory of connections or gauge fields on Hopf-Galois 
  extensions or quantum 
principal bundles were laid in \cite{BrzMaj:gau}.
  Also in  \cite{BrzMaj:gau} the 
first non-trivial example of a connection 
  in a quantum principal bundle, namely 
the q-deformed Dirac magnetic monopole
  potential, was described. 
From
  physics as well as (algebro) geometric point of view it is important to
  understand local structure of quantum principal bundles. A proposal of 
  a 
definition of a locally trivial quantum principal bundle was already
  made in 
\cite{BrzMaj:gau}. These ideas were then developed more fully in 
  the framework 
of $C^*$-algebras by Budzy\'nski and Kondracki in 
  \cite{BudKon:pri} and then 
in a purely algebraic form by Calow and
  Matthes in \cite{CalMat:cov}, 
\cite{CalMat:con}. This has led to the
  formal definition of {\em coverings} 
and {\em complete coverings} of
  non-commutative algebras and also to the 
construction of an explicit
  example of a locally trivial quantum principal bundle 
that was later 
  shown in \cite{HajMat:loc} to be a Hopf-Galois extension. 

In 
  \cite{Brz:str} it has been realised that a general algebraic framework
  for studying Hopf-Galois extensions and their generalisations is provided
  by 
a special class of coalgebras over non-commutative algebras or {\em 
  corings} 
\cite{Swe:pre} termed {\em Galois corings}. On the other hand
  Kontsevich and 
Rosenberg
  \cite{KonRos:smo} have proposed an approach to non-commutative 
algebraic
  geometry, whereby the structure of a non-commutative scheme is 
encoded in
  a coring of a specific kind. Such a coring is termed a {\em cover} 
and a
  non-commutative algebraic space is recovered from a quotient of a 
  derived category of comodules of a cover. Galois corings appear as examples 
  of covers.

The aim of this note is to make clear and explicit
the connection between coverings appearing
  in the context of quantum principal bundles and covers of non-commutative 
  algebraic geometry. The connection is provided by Galois corings. More 
  specifically 
we show that any covering in the sense of Calow and Matthes 
  induces a coring. 
This coring is a Galois coring provided the covering 
  is complete.

\section{Preliminaries}
We work over a field $k$. All algebras
  are over $k$, associative  and with a unit. 
The unadorned
  tensor product is over $k$.

Let $A$ be an algebra. An $A$-bimodule $\cC$ 
  is called an  {\em $A$-coring} provided there exist $A$-bilinear maps $\DC 
  :\cC\to \cC\ot_A\cC$ and  $\eC:\cC\to A$ such that
  $$
(\DC\otimes_A\id_\cC)\circ\DC = (\id_\cC\otimes_A\DC)\circ
\DC, \quad (\eC\otimes_A\id_\cC)\circ\DC =
(\id_\cC\otimes_A\eC)\circ \DC = \id_\cC.
$$
The map $\DC$ is known as a {\em coproduct} and $\eC$ is known as a 
{\em counit} of $\cC$. For a detailed account of the theory of corings the 
reader is referred to \cite{BrzWis:cor}. An example of a coring that is most 
relevant to the subject of this paper can be constructed as follows \cite{Swe:pre}. First, 
consider an algebra map $\iota: B\to A$ (extension of algebras). Then $A$ is 
a $B$-bimodule via $b\cdot a\cdot b' = \iota(b)a\iota(b')$, for all $a\in A$,
$b,b'\in B$. Consequently $A\ot_B A$ is well-defined and an $A$-bimodule with
the actions $a\cdot (a'\ot_B a'')\cdot a''' = aa'\ot_B a''a'''$. Furthermore, 
it is an $A$-coring with coproduct and counit
$$
\Delta_{A\ot_BA}(a\ot_B a') 
= a\ot_B 1_A\ot_A 1_A\ot_B a'= a\ot_B 1_A\ot_B a', \qquad 
\eps_{A\ot_BA}(a\ot_B a') = aa'.
$$
$A\ot_BA$ is known as the canonical or
Sweedler coring associated to a ring extension $\iota:B\to A$.

An element $g$ of an $A$-coring $\cC$ is called a
{\em grouplike element} provided $\DC(g) = g\ot_A g$ and $\eC(g) =1_A$.
Given a grouplike element $g$ one defines
{\em $g$-coinvariants of $A$} by
$$
A^{co\cC}_g= \{b\in A\; |\; b\cdot g = g \cdot
b\}.
$$
$A^{co\cC}_g$ is a subalgebra of $A$, hence the obvious inclusion map gives rise to
the ring extension $A^{co\cC}_g\to A$, and there is the corresponding Sweedler coring
$A\ot_B A$, with $B=A^{co\cC}_g$. A coring $\cC$ with a grouplike element $g$ is called a {\em Galois
coring} provided the $A$-bilinear map
$$
\can_A : A\ot_B A\to \cC, \qquad a\ot_B a'\mapsto a\cdot g\cdot a',
$$
where $B$ are $g$-coinvariants of $A$, is an
isomorphism (of corings).

\section{Complete coverings and Galois corings}
The notion of a {\em complete covering} of a non-commutative algebra was
introduced in \cite{CalMat:cov} as a purely algebraic version of a covering
of a $C^*$-algebra by closed ideals described in the context of quantum
principal bundles in \cite{BudKon:pri}.
\begin{definition}(\cite{CalMat:cov})
Let $B$ be an algebra. A finite family of ideals $(J_i)_{i\in I}$ is called a 
{\em covering} of $B$ if
$$
\bigcap_{i\in I}J_{i}=\{0\}.
$$
\label{def.cover}
\end{definition}
Denote the factor algebras of $B$ with respect to the ideals 
$J_{i}$, $J_{i}+J_{j}$, and so on by $B_{i}$, $B_{ij}$, and  the corresponding canonical surjections by
$$
\pi_{i}:B\to B_{i}, \quad
\pi_{ij}:B\to B_{ij}, \quad
\pi_{ijk}:B\to B_{ijk}, 
$$
etc. Obviously, there exist also surjective maps between various quotients
$$
\pi^{i_1\ldots i_m}_{j_1\ldots j_n}:B_{i_1\ldots i_m}\to B_{i_1\ldots i_mj_1\ldots j_n}, \quad  b+J_{i_1\ldots i_m}\mapsto b+J_{i_1\ldots i_m} + J_{j_1\ldots j_n}.
$$
\begin{definition}(\cite{CalMat:cov}) Given a covering $(J_i)_{i\in I}$  of $B$, define a subalgebra $B_c$ of $\bigoplus_{i\in I}B_{i}$,
$$
B_{c}:=\left.\left\{(a_{i})_{i\in I}\in\bigoplus_{i\in I}B_{i}\right\vert
\pi^{i}_{j}(a_{i})=\pi^{j}_{i}(a_{j})\right\}.
$$
$B_c$ is known as the {\em covering completion} of $B$. A covering $(J_i)_{i\in I}$ is said to be \textit{complete} if the map
$$
\kappa:B\to B_{c}, \qquad \kappa(b)=(\pi_{i}(b))_{i\in I}
$$
is surjective.
\label{def.complete}
\end{definition}

 Note that, since the ideals forming a covering of $B$ intersect trivially, the map $\kappa$ constructed in Definition~\ref{def.complete} is automatically injective. Hence, equivalently, for a covering to be complete one can require $\kappa$ be bijective.

The first main result of this note is that to any  covering 
$(J_i)_{i\in I}$ of an algebra $B$, one can 
associate a coring that captures very closely the way in which 
$B$ is covered by the ideals.

\begin{proposition}\label{prop.cov.cor}
Let $(J_i)_{i\in I}$ be a covering of an algebra $B$ and 
 define an algebra $A=\bigoplus_{i\in I}B_{i}$. 
Consider $\cC=\bigoplus_{i,j}B_{ij}$. 
$\cC$ is an $A$-bimodule with the products 
$$ 
(a_i)_{i\in I}\cdot (a_{jk})_{j,k\in I}\cdot (a'_l)_{l\in I} = 
(\pi^j_k(a_j)a_{jk}\pi^k_j(a'_k))_{j,k\in I}, 
$$
for all $a_i, a'_i\in B_i$ and $a_{jk}\in B_{jk}$. Then $\cC$ is an $A$-coring
with the coproduct 
\begin{eqnarray*}
\DC \big(\big(\pi_{ij}(b_{ij})\big)_{i,j\in I}\big) &=& \sum_{k\in I} \big(\pi_{il}(b_{ik})\big)_{i,l\in I}\ot_A \big(\pi_{mj}(\delta_{kj}1_B)\big)_{m,j\in I}\\
&=&  \sum_{k\in I} \big(\pi_{n}(b_{nk})\big)_{n\in I}\cdot\big(\pi_{il}(1_B)\big)_{i,l\in I}\ot_A \big(\pi_{mj}(\delta_{kj}1_B)\big)_{m,j\in I}\\
&=& \sum_{k\in I} \big(\pi_{il}(\delta_{ik}1_B)\big)_{i,l\in I}\ot_A \big(\pi_{mj}(b_{kj})\big)_{m,j\in I}\\
&=&  \sum_{k\in I} \big(\pi_{il}(\delta_{ik}1_B)\big)_{i,l\in I}\ot_A \big(\pi_{mj}(1_B)\big)_{m,j\in I}\cdot \big(\pi_{n}(b_{kn})\big)_{n\in I},
\end{eqnarray*}
and the counit
$$
\eC \big(\big(\pi_{ij}(b_{ij})\big)_{i,j\in I}\big) = \big(\pi_{i}(b_{ii})\big)_{i\in I}, 
$$
for all $b_{ij}\in B$.
\end{proposition}
\begin{proof} The equivalent forms of the coproduct can easily be verified. It is then clear 
that $\DC$ is an $A$-bimodule map and that $\eC$ is a counit. That $\DC$ is coassociative can also be checked by a straightforward calculation (but, in the light of 
the discussion below this is not necessary).
\end{proof}

There is an algebra inclusion $\iota: B\to A$ given by
$b \mapsto (\pi_{i}(b))_{i\in I}$, i.e.,  of the same explicit form 
as 
the map $\kappa$ in Definition~\ref{def.complete}, hence we can 
consider the corresponding canonical (Sweedler) coring $A\otimes_{B}A$. In fact 
the coring $\cC$ defined in Proposition~\ref{prop.cov.cor} is isomorphic to $A\ot_B A$. The isomorphism can be explicitly defined in two stages as follows.

First, for all $i,j\in I$, there is an $A$-bimodule isomorphism
$$ 
\Phi_{ij}:B_{i}\ot_{B}B_{j}\to B_{ij}, \qquad 
a\ot a'\mapsto\pi^{i}_{j}(a)\pi^{j}_{i}(a'), 
$$ 
with the inverse, for all $b\in B$,
$$
\Phi_{ij}^{-1}: \pi_{ij}(b) \mapsto \pi_i(1_B)\ot_B\pi_j(b) = \pi_i(b)\ot_B \pi_j(1_B).
$$
The $\Phi_{ij}$  taken together give the $A$-bimodule isomorphism 
$$ 
\Phi:\bigoplus_{i,j}B_{i}\ot_B B_{j}\to\bigoplus_{i,j}B_{ij}, 
\qquad \Phi=\bigoplus_{i,j}\Phi_{ij}. 
$$ 
Second, since direct sums commute with tensor products, there is the 
 $A$-bimodule 
isomorphism 
$$ 
\Theta: \bigoplus_{i\in I}B_{i}\ot_{B}\bigoplus_{j\in I}B_{j}\to 
\bigoplus_{i,j\in I}B_{i}\ot_{B}B_{j}, 
\qquad (a_i)_{i\in I}\ot_B (a'_j)_{j\in I}\mapsto (a_i\ot_Ba'_j)_{i,j\in I}. 
$$ 
Combining $\Phi$ with $\Theta$ we obtain the $A$-bimodule isomorphism
\begin{gather*} 
\chi:A\ot_{B}A\to\cC, \qquad \chi=\Phi\circ\Theta; \\ 
(a_{i})_{i\in I}\ot_{B}(a'_{j})_{j\in I}\mapsto 
(\pi^{i}_{j}(a_{i})\pi^{j}_{i}(a'_{j}))_{i,j\in I}, 
\end{gather*} 
with the inverse, for all $b_{ij}\in B$,
$$
\chi^{-1}: \big(\pi_{ij}(b_{ij})\big)_{i,j\in I}\mapsto \sum_{k\in I} \big(\pi_i(b_{ik})\big)_{i\in I}\ot_B \big(\pi_j(\delta_{jk})\big)_{j\in I}.
$$
A straightforward calculation verifies that $\chi$ is an $A$-coring map (and we can deduce from this and the coassociativity of the coproduct in $A\ot_B A$ that the coproduct of $\cC$ is coassociative). 

To understand better how $\cC$ captures the way in which $B$ is covered by $(J_i)_{i\in I}$, we can use  further isomorphisms 
$$A\ot_B A\ot_AA\ot_B A\simeq A\ot_B A\ot_BA,$$
 $$B_i\ot_BB_j\ot_B B_k \to B_{ijk}, \quad a_i\ot_B a_j\ot_B a_k\mapsto  \pi^i_{jk}(a_i)\pi^j_{ik}(a_j)\pi^k_{ij}(a_k)$$
  and 
  $$\bigoplus_{i\in I}B_{i}\ot_{B}\bigoplus_{j\in I}B_{j}\ot_{B}\bigoplus_{k\in I}B_{k}\simeq 
\bigoplus_{i,j,k\in I}B_{i}\ot_{B}B_{j}\ot_B B_k,$$
to identify $\cC\ot_A\cC$ with the $A$-bimodule of ``triple intersections"
$
\cD = \bigoplus_{i,j,k\in I}B_{ijk}.
$
In terms of this identification, the coproduct takes the following simple form, for all $a_{ij}\in B_{ij}$,
$$
\DC \big(\big(a_{ij} \big)_{i,j\in I}\big) = \big(\pi^{ik}_{j}(a_{ik})\big)_{i,j,k\in I}.
$$

As an example of $\cC$ consider first a (topological) space $X$ 
covered by a finite family of
(open) sets $(U_i)_{i\in I}$. Let $B$ be an algebra of functions on $X$, $B = \cO(X)$. The algebra $B$ is then covered by the ideals $J_i = \{f\in   \cO(X)\; |\; \forall x\in X\setminus U_i, \; f(x) =0\}$, so that
$B_i = \cO(U_i)$, and each of the $\pi_i$ is simply 
the restriction of a function on
$X$ to a function on $U_i$, $\pi_i(f) = f\mid_{U_i}$. In addition, $B_{ij} = \cO(U_i\cap U_j)$. The algebra $A = \bigoplus_{i\in I} B_i$ can be identified with the algebra of functions on the disjoint union $\bigsqcup_{i\in I} U_i$, $A=\cO\big(\bigsqcup_{i\in I} U_i\big)$. Furthermore, 
$\cC = \bigoplus_{i,j\in I} B_{ij}$ can be identified with the algebra of functions on the disjoint union of intersections $\bigsqcup_{i,j\in I} U_i\cap U_j$, i.e., $\cC=\cO\big(\bigsqcup_{i,j\in I} U_i\cap U_j\big)$, and also $\cD = \cO\big(\bigsqcup_{i,j,k\in I} U_i\cap U_j\cap U_k\big)$. The coproduct and the counit come out as, for all $f\in \cO\big(\bigsqcup_{i,j\in I} U_i\cap U_j\big)$,
$$
\DC(f)\mid_{U_i \cap U_j \cap U_k} = f_{ik}\mid_{U_i\cap U_j \cap U_k}, \qquad \eC(f)\mid_{U_i} = f\mid_{U_i\cap U_i},
$$
where $f_{ik} = f \mid_{U_i\cap U_k}$. In this way we obtain an example of  coring described in \cite{KonRos:smo}.

The following theorem establishes a relationship between 
 the notions of complete coverings and Galois corings and also gives the criterion, when a covering is complete.
\begin{theorem} \label{thm.main}
Let  $(J_{i})_{i\in I}$ be a covering of an algebra $B$, $B_i = B/J_i$, and let $A=\bigoplus_{i\in I} B_i$. Let  $\cC=\bigoplus_{i,j}B_{ij}$ be the associated $A$-coring constructed in Proposition~\ref{prop.cov.cor}. Set $g = \big(\pi_{ij}(1_B)\big)_{i,j\in I}$. 
\begin{zlist}
\item If $(J_{i})_{i\in I}$  is a complete covering of  $B$, then $(\cC,g)$ is a Galois coring.
\item If $A$ is a faithfully flat left or right $B$ module, then $(J_{i})_{i\in I}$  is a complete covering of  $B$.
\end{zlist}
\end{theorem} 
\begin{proof} 
First note that $g$ is a grouplike element and that it is 
the image of the grouplike element $1_{A}\ot_{B}1_{A}$ of 
the Sweedler $A$-coring $A\ot_B A$ under the isomorphism $\chi$. Since $\pi_{ij}(1_B)$ is the unit in $B_{ij}$, the forms of the right and left multiplications in $\cC$ immediately imply that, for all $(a_i)_{i\in I}\in A$,
$$
(a_i)_{i\in I}\cdot g = (a_i)_{i\in I}\cdot (\pi_{jk}(1_B))_{j,k\in I} = (\pi^j_k(a_j))_{j,k\in I}, 
$$
and
$$
 g\cdot (a_i)_{i\in I} = (\pi_{jk}(1_B))_{j,k\in I}\cdot (a_i)_{i\in I} = (\pi^k_j(a_k))_{j,k\in I}.
$$
Therefore $(a_i)_{i\in I}\in A^{co\cC}_g$ if and only if $(a_i)_{i\in I}\in B_c$, where $B_c$ is a covering completion of $B$. Thus 
$$
B_c = A^{co\cC}_g.
$$
In view of the isomorphism $\chi$ this identification of $B_c$ as $g$-coinvariants  immediately implies assertion (1).

The assertion (2) follows by \cite[Theorem~4.12]{CaeDeG:com} (cf.\ \cite[Theorem~3.10]{KaoGom:com}). We present a 
direct proof based on the proof of  \cite[Lemma~4.2]{Tak:ext}. Recall that $A^{co\cC}_g = B_c$ and that $B$ is a subalgebra of $B_c$ via $\kappa: B\to B_c$. Note that
$$
\chi(A\ot_B B_c) = A\cdot g, \qquad \chi(B_c\ot_B A) = g\cdot A.
$$
Note further that $g\cdot A\simeq A$ as a $(B,A)$-bimodule and $A\cdot g \simeq A$ as an $(A,B)$-bimodule via the counit $\eC$. Since $\chi$ is an isomorphism, the map
$$
\xymatrix{
A\ot_B B_c \ar[rr]^{\chi} & & A\cdot g \ar[rr]^{\eC} && A},\qquad a\ot_B b\mapsto ab
$$ 
is an isomorphism of $(A,B)$-bimodules, and the map
$$
\xymatrix{
B_c\ot_B A \ar[rr]^{\chi} & & g\cdot A \ar[rr]^{\eC} && A},\qquad b\ot_B a\mapsto ba
$$ 
is an isomorphism of $(B,A)$-bimodules. These lead to the factorisations of 
the identity maps
$$
\xymatrix{A \ar[r]^{\simeq} & B\ot_B A \ar[r]^{\kappa\ot_B A} & B_c\ot_B A \ar[r]^{\simeq} & A}, \quad \xymatrix{A \ar[r]^{\simeq} & A\ot_B B \ar[r]^{A\ot_B\kappa} & A\ot_B B_c \ar[r]^{\simeq} & A,}
$$
implying that $\kappa\ot_B A$ and $A\ot_B\kappa$ are isomorphisms. Thus if $A$ is faithfully flat as either left or right $B$-module, then $\kappa: B\to B_c$ is an isomorphism, hence $(J_{i})_{i\in I}$  is a complete covering of  $B$, as claimed.
\end{proof}
\begin{corollary}\label{cor1}
Let  $(J_{i})_{i\in I}$ be a complete covering of an algebra $B$, $B_i = B/J_i$, and let $A=\bigoplus_{i\in I} B_i$. If $A$ is a projective left (resp.\ right) $B$-module, then $A$ is a faithfully flat left (resp.\ right) $B$-module.
\end{corollary}
\begin{proof} By Theorem~\ref{thm.main}(1), $\cC=\bigoplus_{i,j\in I} B_{ij}$ is a 
Galois $A$-coring, hence $A$ is a principal right (resp.\ left) $\cC$-comodule in the sense of \cite[Definition~4.1]{Brz:gal}. By \cite[Theorem~4.3]{Brz:gal} (resp.\ the left-handed version of \cite[Theorem~4.3]{Brz:gal}), $A$ is a faithfully flat left (resp.\ right) $B$-module.
\end{proof}

In particular, Theorem~\ref{thm.main} implies that if  $(J_{i})_{i\in I}$ is a covering of an algebra $B$ such that $A$ is a faithfully flat left $B$-module (e.g., if $(J_{i})_{i\in I}$ is a complete covering and ${}_BA$ is projective by Corollary~\ref{cor1}), then $\cC$ is a flat left $A$-module and the induction functor $-\ot_B A$ is an equivalence from the category of right $B$-modules to the category of right $\cC$-comodules by \cite[Theorem~5.6]{Brz:str}. Following Kontsevich and Rosenberg, $\cC$-comodules are understood as quasi-coherent sheaves on a non-commutative space. Since the category of comodules (quasi-coherent sheaves) is equivalent to a module category, the underlying space is a  non-commutative affine space with the 
coordinate ring $B$ (cf.\ \cite[p.\ 2134]{Smi:sub}).

\end{document}